\documentclass[11pt,reqno]{amsart}
\usepackage[T1]{fontenc}
\usepackage{lmodern}
\usepackage[letterpaper,margin=1in]{geometry}
\usepackage{amsmath,amssymb,amsthm,mathtools}
\usepackage{microtype}
\usepackage[colorlinks=true,linkcolor=black,citecolor=black,urlcolor=blue]{hyperref}
\allowdisplaybreaks[2]
\newcommand{\R}{\mathbb R}
\newcommand{\B}{\mathbb B^d}
\newcommand{\Pn}{\mathcal P_n}
\newcommand{\norm}[1]{\lVert#1\rVert}
\newcommand{\om}{\omega_2}
\newcommand{\T}{\mathcal T}
\newcommand{\op}{\mathrm{op}}
\newtheorem{theorem}{Theorem}[section]
\newtheorem{lemma}[theorem]{Lemma}
\newtheorem{corollary}[theorem]{Corollary}

\title[Convex polynomial approximation on the unit ball]{Convex polynomial approximation on the unit ball}
\author{Andriy Bondarenko}
\address{Department of Mathematical Sciences, Norwegian University of
Science and Technology, Trondheim, Norway}
\email{andriybond@gmail.com}
\thanks{Andriy Bondarenko acknowledges partial support from the Research
Council of Norway (Grant~334466).}
\author{Andriy Prymak}
\address{Department of Mathematics, University of Manitoba,
Winnipeg, MB R3T~2N2, Canada}
\email{prymak@gmail.com}
\thanks{Andriy Prymak acknowledges support from the Natural Sciences and
Engineering Research Council of Canada (Discovery Grant~RGPIN-2026-06488).}
\author{Igor Shevchuk}
\address{Faculty of Mechanics and Mathematics, Taras Shevchenko National University of Kyiv, 01601 Kyiv, Ukraine}
\email{shevchukio@knu.ua}
\thanks{Igor Shevchuk acknowledges support from Simons Foundation
(\mbox{PD-COE-00031486}).}
\date{}
\subjclass[2020]{41A10, 41A25, 41A29}
\keywords{convex approximation, algebraic polynomials, second modulus of smoothness, Jackson-type estimate, unit ball}

\begin{document}
\begin{abstract}
Let $\B$ be the closed unit ball in $\R^d$, $d\ge1$.
We prove that, for every convex $f\in C(\B)$ and every integer $n\ge1$,
there exists a polynomial $P_n$ of total degree at most $n$, convex
on all of $\R^d$, such that
$\norm{f-P_n}_{C(\B)}\le c\,\om(f,n^{-1};\B)$.
Here $\om$ is the ordinary second modulus of smoothness and $c=c_d$
depends only on $d$. In particular, we obtain the rate $O(n^{-2})$ for
convex functions with uniformly bounded second derivatives.
\end{abstract}
\maketitle

\section{Introduction and main result}

Shvedov~\cite{Shvedov1981} proved that any continuous convex function on
an interval admits convex polynomial approximants with a Jackson
estimate in the uniform norm in terms of the second modulus of
smoothness. The univariate theory has since seen substantial progress,
including pointwise estimates and estimates involving higher-order
moduli of smoothness. We refer to the survey
\cite{KopotunLeviatanPrymakShevchuk2011} for further details.

The multivariate problem is substantially different. Shvedov
\cite{Shvedov1982} established a first-modulus estimate for convex
polynomial approximation on convex bodies. The validity of a
second-modulus bound in several variables remained open. Related
results in multivariate convexity-preserving approximation can be found
in \cite{BondarenkoPrymak2004,Budnik1994,KairePrymak2025,KonovalovKopotunMaiorov2010,MaizlishPrymak2015}.

We prove an estimate involving the ordinary second modulus for all
continuous convex functions on the Euclidean unit ball in every fixed
dimension, with approximating polynomials convex on the whole space.

Fix an integer $d\ge1$. Here and throughout, $|x|$ is the Euclidean norm. Put
$\B=\{x\in\R^d:|x|\le1\}$, write
$B(0,1)=\{x\in\R^d:|x|<1\}$ for the open unit ball, and let $\Pn$
denote the polynomials of total degree at most $n$. For $f\in C(\B)$, write
$\norm{f}_{\B}=\max_{x\in\B}|f(x)|$ and
\begin{equation}\label{eq:modulus}
 \om(f,t;\B)=\sup_{\substack{|u|\le t\\x-u,x,x+u\in\B}}
 \bigl|f(x+u)-2f(x)+f(x-u)\bigr|.
\end{equation}
For a real $d\times d$ matrix $A$, its Euclidean operator norm is
$\norm A_{\op}=\sup_{|v|=1}|Av|$.
By $C^2(\B)$ we mean functions that are twice continuously differentiable
in the interior of $\B$, with all partial derivatives of order at most
two extending continuously to the closed ball. For $h\in C^2(\B)$, set
\[
 \norm{D^2h}_{\B,\op}
 =\sup_{x\in\B}\norm{D^2h(x)}_{\op}
 =\sup_{\substack{x\in\B\\|v|=1}}|v^\top D^2h(x)v|,
\]
where the last equality uses symmetry of the Hessian $D^2h(x)$.
Throughout the paper, $c>0$ denotes a constant depending only on $d$,
not necessarily the same at different occurrences; $a\asymp b$ means
$cb\le a\le c^{-1}b$ for some such $c>0$. Constants called absolute
are independent of $d$.

\begin{theorem}\label{thm:main}
For every convex $f\in C(\B)$ and every integer $n\ge1$, there is
$P_n\in\Pn$, convex on $\R^d$, such that
\begin{equation}\label{eq:main}
 \norm{f-P_n}_{\B}\le c\,\om(f,n^{-1};\B).
\end{equation}
\end{theorem}

\begin{corollary}\label{cor:hessian}
If $f\in C^2(\B)$ is convex, there are polynomials $P_n\in\Pn$,
convex on $\R^d$, for which
\[
 \norm{f-P_n}_{\B}\le cn^{-2}\norm{D^2f}_{\B,\op}.
\]
In particular, if $|\partial_i\partial_j f|\le1$ on $\B$ for all
$1\le i,j\le d$, then
$\norm{f-P_n}_{\B}\le cn^{-2}$.
\end{corollary}

Indeed, $\om(f,t;\B)\le t^2\norm{D^2f}_{\B,\op}$, and the stated
entrywise bound gives $\norm{D^2f}_{\B,\op}\le d$.
Our proof uses a positive convolution on the ball, related to
\cite{Xu2005}, with a polynomial kernel constructed explicitly below.
A radial correction, in the
spirit of \cite{MaizlishPrymak2015}, makes the approximant globally
convex without changing the order of approximation.

\section{Two approximation lemmas}

The estimates in the following lemma, up to constants depending only
on $d$, follow from Dekel's theorem on the equivalence of the modulus
of smoothness and the Sobolev $K$-functional over convex domains
\cite[Theorem~1.1]{Dekel2010}. We include a short direct proof on the
ball, valid up to the boundary, that gives the explicit approximation
constant $5/2$.

\begin{lemma}\label{lem:smoothing}
For each $f\in C(\B)$ and $0<t\le\tfrac14$, there is a function $h_t$,
of class $C^\infty$ on an open neighborhood of $\B$, satisfying
\begin{equation}\label{eq:smoothing}
 \norm{f-h_t}_{\B}\le\tfrac52\om(f,t;\B),\qquad
 \norm{D^2h_t}_{\B,\op}\le ct^{-2}\om(f,t;\B).
\end{equation}
\end{lemma}
\begin{proof}
Choose an even nonnegative $\phi\in C_c^\infty(B(0,1))$ of integral
one, put $\varepsilon=t/2$ and
$\phi_\varepsilon(u)=\varepsilon^{-d}\phi(u/\varepsilon)$, and define
\[
 \begin{aligned}
 g_t(z)&=\int\phi_\varepsilon(u)f(z+u)\,du,
 \qquad |z|<1-\varepsilon,\\
 h_t(x)&=2g_t((1-t)x)-g_t((1-2t)x).
 \end{aligned}
\]
The latter is of class $C^\infty$ on an open neighborhood of $\B$. Evenness gives
$|g_t(z)-f(z)|\le\frac12\om(f,\varepsilon;\B)$, while
$x,(1-t)x,(1-2t)x$ are equally spaced in $\B$ with step at most $t$.
This proves the first estimate in \eqref{eq:smoothing}.
For each unit vector $v$, we have
\[
 D_v^2g_t(z)=\frac12\int D_v^2\phi_\varepsilon(u)
 \bigl(f(z+u)-2f(z)+f(z-u)\bigr)\,du.
\]
Since $\norm{D_v^2\phi_\varepsilon}_{L^1}\le c\varepsilon^{-2}$,
the chain rule proves the second estimate.
\end{proof}

The smoothing function $h_t$ need not be convex.

For the kernel construction and the global correction, fix
\[
 q_d=\left\lceil\frac d2\right\rceil+4,
 \qquad 2q_d>d+7.
\]
For $d=2$, this gives $q_d=5$ and the kernel power $10$.
The crucial result is given by the next lemma.

\begin{lemma}\label{lem:operator}
For every integer $m\ge1$ there is a positive linear operator
$\T_m:C(\B)\to\mathcal P_{2q_d m+2}$, preserving constants and convexity
on $\B$, such that for every function $h$ of class $C^2$ on an open
neighborhood of $\B$
\begin{align}
 \norm{D^2\T_mh}_{\B,\op}&\le\norm{D^2h}_{\B,\op},
 \label{eq:hesscontract}\\
 \norm{(I-\T_m)(h-\ell_h)}_{\B}&\le
 cm^{-2}\norm{D^2h}_{\B,\op},
 \quad\ell_h(x)=h(0)+\nabla h(0)\cdot x.
 \label{eq:operatorerror}
\end{align}
\end{lemma}

\noindent\emph{Idea of the proof.}
We construct the operator by averaging over $S^{d+1}$, the unit sphere
in $\R^{d+2}$, whose projection onto the first $d$ coordinates is the ball.
The kernel will concentrate the integral near the point $x$ being
approximated. More precisely, we will estimate its vector first moment
and scalar second moment about $x$.
To preserve convexity, we choose the kernel as a double integral of a
nonnegative polynomial. An integration-by-parts identity then transfers
second derivatives of the approximant to second derivatives of the
original function, through a positive operator. The moment estimates
and Taylor's formula give the approximation estimate after subtracting
an affine function. Finally, the action of the operator on $|x|^2$
gives the Hessian contraction estimate.

\begin{proof}
\emph{Step 1: the spherical lift and its measure.}
Let $V_d=\pi^{d/2}/\Gamma(d/2+1)$ be the volume of $\B$, and let
$w(x)=1-|x|^2$. Define the lift $\tau:\B\to S^{d+1}$ by
\[
 \tau(x)=(x,\sqrt{w(x)},0),\qquad x\in\B.
\]
To use $h\in C(\B)$ in an integral over the sphere, put
\[
 \widetilde h(y)=h(y_1,\ldots,y_d),\qquad y\in S^{d+1}.
\]
Thus $\widetilde h$ takes the value of $h$ at the projection of $y$
onto the ball, independently of the last two coordinates. In particular,
$\widetilde h(\tau(x))=h(x)$.
Let $dS$ be surface area measure on $S^{d+1}$, and put
\[
 d\sigma(y)=\frac{dS(y)}{2\pi V_d},\qquad
 \int_{S^{d+1}}d\sigma(y)=1.
\]
For a univariate polynomial kernel $K$, define
\begin{equation}\label{eq:sphere}
 \T_{0,K}h(x)=\int_{S^{d+1}}K(\tau(x)\cdot y)
                   \widetilde h(y)\,d\sigma(y).
\end{equation}

To express this operator on the ball, parametrize almost every point of
$S^{d+1}$ by
\[
 y=(z,\sqrt{w(z)}\cos\phi,\sqrt{w(z)}\sin\phi),
 \qquad |z|<1,\quad 0\le\phi<2\pi.
\]
Here $z=(y_1,\ldots,y_d)$ is the projection of $y$ onto the ball.
The Gram matrix of this parametrization has $z$-block
$I_d+zz^\top/w(z)$, angular entry $w(z)$, and zero mixed entries.
Its determinant is one, so the surface element is $dz\,d\phi$ and
\begin{equation}\label{eq:sphere-measure}
 d\sigma(y)=\frac{dz}{V_d}\frac{d\phi}{2\pi}.
\end{equation}
Thus the projected measure is normalized volume on $\B$. Moreover,
\[
 \tau(x)\cdot y=x\cdot z+\sqrt{w(x)w(z)}\cos\phi.
\]
Integrating in $\phi$ and substituting $u=\cos\phi$ give
\[
 \T_{0,K}h(x)=\frac1{V_d}\int_{\B} h(z)
 \left[\frac1\pi\int_{-1}^1
 \frac{K\bigl(x\cdot z+\sqrt{w(x)w(z)}\,u\bigr)}{\sqrt{1-u^2}}
 \,du\right]dz.
\]
The inner integral is obtained by averaging over the circle above $z$.
It is the case $\alpha=0$ of the family introduced next.

\smallskip
\noindent\emph{Step 2: weighted operators and convexity preservation.}
Differentiation shifts the ball weight, so we introduce an auxiliary
parameter $\alpha\ge0$. For a univariate polynomial $K$ and
$h\in C(\B)$, define
\begin{align*}
 c_\alpha&=\frac{\Gamma(\alpha+1)}{\sqrt\pi\,\Gamma(\alpha+1/2)},\\
 A_{\alpha,K}(x,z)&=c_\alpha\int_{-1}^1
 K\bigl(x\cdot z+\sqrt{w(x)w(z)}\,u\bigr)
 (1-u^2)^{\alpha-1/2}\,du,\\
 b_{d,\alpha}&=\frac{\Gamma(\alpha+d/2+1)}{\pi^{d/2}\Gamma(\alpha+1)},\\
 \T_{\alpha,K}h(x)&=b_{d,\alpha}
 \int_{\B} h(z)w(z)^\alpha A_{\alpha,K}(x,z)\,dz.
\end{align*}
Both $c_\alpha(1-u^2)^{\alpha-1/2}\,du$ and
$b_{d,\alpha}w(z)^\alpha\,dz$ are nonnegative measures of total mass one.
For $\alpha=0$, these formulas agree with the spherical integral in
\eqref{eq:sphere}. Only the values $\alpha=0,1,2$ will be needed.
Odd powers of $u$ vanish in the inner integral, so the square roots
occur only to even powers. Consequently,
$\T_{\alpha,K}h\in\mathcal P_{\deg K}$.
Also, $\T_{\alpha,K}$ is positive whenever $K\ge0$ on $[-1,1]$.

We now transfer derivatives from the kernel to the function being
averaged. Assume that $h$ is of class $C^2$ on an open neighborhood
of $\B$. For $x,z$ in the interior of $\B$, differentiation and
integration by parts in $u$ give
\[
 w(z)^\alpha\partial_{x_i}A_{\alpha,K}(x,z)
 =-\frac1{2(\alpha+1)}\partial_{z_i}
 \bigl(w(z)^{\alpha+1}A_{\alpha+1,K'}(x,z)\bigr).
\]
The boundary factor vanishes. Integration by parts in $z$ then yields
the first identity below; applying it twice yields the second:
\begin{equation}\label{eq:diff}
 \begin{aligned}
 \partial_i\T_{\alpha,K}h
 &=\frac1{2\alpha+d+2}\T_{\alpha+1,K'}(\partial_i h),\\
 D_v^2\T_{0,K}h&=\frac1{(d+2)(d+4)}\T_{2,K''}(D_v^2h).
 \end{aligned}
\end{equation}
Both sides are polynomials in $x$, so the identities hold on the
boundary as well. A related weight-shifting differentiation identity for
orthogonal polynomial partial sums is given in
\cite[Lemma~3.1, equation~(3.3)]{PinarXu2018}.
For our kernel operators, \eqref{eq:diff} has been established directly,
including the normalization factor.

The second identity explains our condition on $K$: if $K''\ge0$ on
$[-1,1]$, then positivity of $\T_{2,K''}$ implies that
$D_v^2\T_{0,K}h\ge0$ whenever $D_v^2h\ge0$. Hence $\T_{0,K}$ preserves
convexity. The conclusion extends to continuous convex functions by
inward contraction, mollification, and uniform convergence.

\smallskip
\noindent\emph{Step 3: choosing and normalizing the kernel.}
We construct the kernel and verify all the required properties directly.
Let $T_j(s)=\cos(j\arccos s)$, $-1\le s\le1$, be the Chebyshev
polynomial of degree $j$, and set
\begin{equation}\label{eq:kernel}
 R_m(s)=1+2\sum_{j=1}^mT_j(s),\qquad
 J_m(t)=\int_{-1}^t(t-s)R_m(s)^{2q_d}\,ds.
\end{equation}
Put
\begin{gather*}
 a_d=\frac{\Gamma((d+2)/2)}{\sqrt\pi\,\Gamma((d+1)/2)},\\
 Z_m=a_d\int_{-1}^1J_m(t)(1-t^2)^{(d-1)/2}\,dt,
 \qquad K_m=J_m/Z_m.
\end{gather*}
This construction is due to Cao and Gonska
\cite[formulas~(1.4) and~(4.1)]{CaoGonska1994}: up to a positive
scalar factor, $J_m$ is their integrated higher-order Jackson kernel
$F_{2m+1,q_d,2}$. The estimates needed for the spherical weight are
proved below. Alternatively, it is possible to adopt the generating-polynomial construction of Kopotun, Leviatan, Prymak
and Shevchuk
\cite[Section~2.4, Lemma~2.13]{KopotunLeviatanPrymakShevchuk2016}
by modifying the power and normalization for the spherical
weight.

\begin{samepage}
Double integration gives $K_m,K_m''\ge0$ on $[-1,1]$, while
$\deg K_m\le2q_d m+2$. The normalization gives, for all
$x\in\B$,
\begin{equation}\label{eq:kernel-mass}
 \int_{S^{d+1}}K_m(\tau(x)\cdot y)\,d\sigma(y)
 =a_d\int_{-1}^1K_m(t)(1-t^2)^{(d-1)/2}\,dt=1.
\end{equation}
\end{samepage}
The first equality is the spherical integration formula for functions
of a single coordinate. It follows from Steps~1--2 that
$\T_m:=\T_{0,K_m}$ is positive, preserves constants and convexity, and
has the required degree bound.

We record the concentration estimates with this dimension-dependent
spherical weight. The identity
\[
 R_m(\cos\theta)
 =\frac{\sin((m+\tfrac12)\theta)}{\sin(\theta/2)}
\]
gives
\[
 |R_m(1-u)|\le c\min(m,u^{-1/2}),\qquad
 R_m(1-u)\ge cm\quad(0\le u\le c_0m^{-2}).
\]
Here $c_0,c_1>0$ are sufficiently small fixed constants.
The upper bound holds for $0<u\le2$. Writing
$J_m(1-u)=\int_u^2(v-u)R_m(1-v)^{2q_d}\,dv$ therefore yields
\begin{align*}
 J_m(1-u)&\le c\min(m^{2q_d-4},u^{2-q_d}),&&0<u\le2,\\
 J_m(1-u)&\ge cm^{2q_d-4},&&0\le u\le c_1m^{-2}.
\end{align*}
For the lower bound one may restrict the defining integral to an interval
of length comparable to $m^{-2}$ within the region where $R_m\ge cm$.
Integrating these bounds against
$u^{(d-1)/2}(2-u)^{(d-1)/2}$ gives
\begin{equation}\label{eq:momentsK}
 \begin{split}
 Z_m&\asymp m^{2q_d-d-5},\\
 \delta_m&:=a_d\int_{-1}^1(1-t)K_m(t)(1-t^2)^{(d-1)/2}\,dt
          \le cm^{-2}.
 \end{split}
\end{equation}
Indeed, replacing $K_m$ by $J_m$ in the integral defining $\delta_m$
and splitting at $u=m^{-2}$ gives a bound $cm^{2q_d-d-7}$.
The inequality $2q_d>d+7$ ensures this bound for the tail integral.
Dividing by $Z_m$ gives $\delta_m\le cm^{-2}$, which explains the
dimension-dependent choice of the kernel power.

\smallskip
\noindent\emph{Step 4: moment estimates and the approximation error.}
Fix $x\in\B$. The ball representation from Steps~1--2 gives
\[
 \T_mh(x)=\int_{\B} h(z)\,d\nu_x(z),\qquad
 d\nu_x(z)=\frac1{V_d} A_{0,K_m}(x,z)\,dz.
\]
The measure $\nu_x$ is nonnegative and has total mass one by
\eqref{eq:kernel-mass}.
Rotational symmetry and the definition of $\delta_m$ give
\[
 \int_{S^{d+1}}yK_m(\tau(x)\cdot y)\,d\sigma(y)
 =(1-\delta_m)\tau(x).
\]
Also, since $|y-\tau(x)|^2=2(1-\tau(x)\cdot y)$ on $S^{d+1}$,
\[
 \int_{S^{d+1}}|y-\tau(x)|^2K_m(\tau(x)\cdot y)\,d\sigma(y)
 =2\delta_m.
\]
Projecting onto the first $d$ coordinates therefore yields
\begin{equation}\label{eq:moments}
 \int_{\B}(z-x)\,d\nu_x(z)=-\delta_m x,\qquad
 \int_{\B}|z-x|^2\,d\nu_x(z)\le2\delta_m.
\end{equation}
These are the first vector moment and the second scalar moment of the
kernel about the evaluation point $x$. The norm of the first moment and
the second moment are both $O(m^{-2})$, uniformly in $x$.

To obtain \eqref{eq:operatorerror}, put $g=h-\ell_h$ and
$M=\norm{D^2h}_{\B,\op}$. Then $\nabla g(0)=0$, so
$|\nabla g(x)|\le M$ on $\B$. Taylor's formula gives
\[
 |g(z)-g(x)-\nabla g(x)\cdot(z-x)|\le\tfrac12M|z-x|^2.
\]
Integrating with respect to $\nu_x$ and using \eqref{eq:moments}, we obtain
\begin{align*}
 |\T_mg(x)-g(x)|
 &\le|\nabla g(x)|\left|\int_{\B}(z-x)\,d\nu_x(z)\right|\\
 &\quad+\tfrac12M\int_{\B}|z-x|^2\,d\nu_x(z)\\
 &\le2M\delta_m.
\end{align*}
Together with \eqref{eq:momentsK}, this proves the approximation estimate.

\smallskip
\noindent\emph{Step 5: the Hessian contraction.}
We finish by checking how $\T_m$ acts on the quadratic function $|x|^2$.
Let
\[
 \eta_m=\int_{S^{d+1}}(\tau(x)\cdot y)^2K_m(\tau(x)\cdot y)\,d\sigma(y),
\]
which is independent of $x$ by rotational symmetry. Symmetry around
the axis through $\tau(x)$ gives
\begin{align*}
 \int_{S^{d+1}}yy^\top K_m(\tau(x)\cdot y)\,d\sigma(y)
 &=\eta_m\tau(x)\tau(x)^\top\\
 &\quad+\frac{1-\eta_m}{d+1}\bigl(I_{d+2}-\tau(x)\tau(x)^\top\bigr),
\end{align*}
where $I_{d+2}$ is the $(d+2)\times(d+2)$ identity matrix. Adding the first $d$
diagonal entries yields
\[
 \T_m(|\cdot|^2)(x)=\frac{d(1-\eta_m)}{d+1}+b_m|x|^2,
 \qquad b_m=\frac{(d+2)\eta_m-1}{d+1}.
\]
Convexity preservation gives $b_m\ge0$, and \eqref{eq:kernel-mass} gives
$\eta_m\le1$, hence $b_m\le1$.
The functions $M|x|^2/2\pm h$ are convex on $\B$.
Applying $\T_m$ and differentiating their images twice gives
\[
 b_mM|v|^2\pm D_v^2\T_mh(x)\ge0,
 \qquad x\in\B,\quad v\in\R^d.
\]
Thus $|D_v^2\T_mh(x)|\le M|v|^2$, which is
\eqref{eq:hesscontract} and completes the proof.
\end{proof}

\section{A global convexity correction}

The strategy of making a locally convex polynomial globally convex by
adding a radial polynomial is classical; see
\cite[proof of Theorem~2]{Shvedov1982} and
\cite{MaizlishPrymak2015}. Shvedov's construction proves density by
choosing the correction degree in terms of the approximant and the
required error. For a related quadratic correction that restores
convexity on a fixed cube, see
\cite[Section~2.2.1]{KonovalovKopotunMaiorov2010}.
Here we need simultaneous control of the correction degree, its size on
$\B$, and its radial and tangential second derivatives outside $\B$.
We therefore give the following elementary estimates, which will yield
a correction of degree $O(m)$ and an error bounded by the right-hand
side of \eqref{eq:Kfunctional}.

We shall use $T_N(r)=\cosh(N\operatorname{arcosh}r)$ for $r\ge1$.
For $d\ge2$, given a unit vector $e$ and a vector $v$, choose a unit vector
$e^\perp\perp e$ so that $v=ae+be^\perp$. If $v$ is parallel to $e$,
put $b=0$. For $d=1$, only the radial term $a$ is needed below.

\begin{lemma}\label{lem:outside}
If $p\in\mathcal P_N$, $N\ge2$, and $\norm p_{\B}\le\varepsilon$,
then, for $r\ge1$,
\begin{equation}\label{eq:outside}
 |D_v^2p(re)|\le c\varepsilon T_N(r)
 \bigl(N^4a^2+N^2b^2\bigr).
\end{equation}
\end{lemma}
\begin{proof}
We use the classical Markov inequality on $[-1,1]$, the
trigonometric Bernstein inequality, and the classical exterior estimate
\begin{equation}\label{eq:Cheb}
 |u(r)|\le\norm u_{[-1,1]}T_N(r),\qquad
 u\in\mathcal P_N(\R),\quad r\ge1.
\end{equation}

For $d\ge2$, choose orthonormal coordinates in $\R^d$ with
$e=e_1$ and $e^\perp=e_2$, and put
$\Omega=x_1\partial_2-x_2\partial_1$.
Bernstein's inequality on circles inside $\B$ gives
$\norm{\Omega p}_{\B}\le N\varepsilon$ and
$\norm{\Omega^2p}_{\B}\le N^2\varepsilon$.
For $g(s)=p(se_1)$ and $k(s)=(\Omega p)(se_1)$, Markov's
inequality followed by \eqref{eq:Cheb} gives
\begin{align*}
 |g'(r)|&\le N^2\varepsilon T_N(r),&
 |g''(r)|&\le N^4\varepsilon T_N(r),\\
 |k(r)|&\le N\varepsilon T_N(r),&
 |k'(r)|&\le N^3\varepsilon T_N(r).
\end{align*}
Also $|\Omega^2p(re_1)|\le N^2\varepsilon T_N(r)$.
At $re_1$,
\[
 p_{11}=g'',\qquad p_{12}=\frac{k'}r-\frac{k}{r^2},
 \qquad p_{22}=\frac{\Omega^2p}{r^2}+\frac{g'}r.
\]
These identities and $2N^3|ab|\le N^4a^2+N^2b^2$ prove
\eqref{eq:outside}, with an absolute constant. For $d=1$, Markov's
inequality and \eqref{eq:Cheb} directly give the radial estimate.
\end{proof}

\begin{lemma}\label{lem:radial}
For each integer $l\ge1$, there is a polynomial $Q_l$ of degree
$4q_d l+2$ such that $Q_l',Q_l''\ge0$ on $[0,\infty)$,
\begin{equation}\label{eq:Qsmall}
 0\le Q_l(s)\le cl^{-2},\qquad 0\le s\le1,
\end{equation}
and, for $r\ge1$,
\begin{equation}\label{eq:Qlarge}
 Q_l'(r^2)\ge c^{-1}T_{2q_d l}(r),\qquad
 Q_l''(r^2)\ge c^{-1}l^2T_{2q_d l}(r).
\end{equation}
\end{lemma}
\begin{proof}
Using $R_l$ from \eqref{eq:kernel}, define
\[
 B_l=\int_0^1R_l(2s-1)^{4q_d}\,ds,\qquad
 Q_l(s)=B_l^{-1}\int_0^s(s-u)R_l(2u-1)^{4q_d}\,du.
\]
The elementary identity
$R_l(\cos\theta)=\sin((l+\tfrac12)\theta)/\sin(\theta/2)$ gives
\begin{align*}
 |R_l(1-u)|&\le c\min(l,u^{-1/2}),&&0<u\le2,\\
 R_l(1-u)&\ge cl,&&0\le u\le c_0l^{-2}.
\end{align*}
Consequently, $B_l\asymp l^{4q_d-2}$ and
$\int_0^1(1-u)R_l(2u-1)^{4q_d}\,du\le cl^{4q_d-4}$.
This proves \eqref{eq:Qsmall}; also $Q_l'(1)=1$.

Put $u=\operatorname{arcosh}r\ge0$. Since
$R_l(\cosh 2u)\ge l\cosh(lu)$, we have
\[
 Q_l''(\cosh^2u)\ge c^{-1}l^2\cosh^{4q_d}(lu)
 \ge c^{-1}l^2\cosh(2q_d lu).
\]
If $u\le1/l$, then $Q_l'(\cosh^2u)\ge1$ and
$\cosh(2q_d lu)\le\cosh(2q_d)$.
If $u\ge1/l$, restrict the integral defining $Q_l'(\cosh^2u)$
to $s=\cosh^2v$, $u/2\le v\le u$. Since
$ds=\sinh(2v)\,dv\ge2v\,dv$,
\[
 Q_l'(\cosh^2u)\ge c^{-1}l^2u^2\cosh^{4q_d}(lu/2)
 \ge c^{-1}\cosh(2q_d lu).
\]
Here we used $\cosh^{4q_d}z\ge2^{1-4q_d}\cosh(4q_d z)$.
\end{proof}

Put $L_d=10q_d+1$.

\begin{lemma}\label{lem:global}
Let $f\in C(\B)$ be convex and let $h$ be of class $C^2$ on an open
neighborhood of $\B$.
For $m\ge4$, there is a polynomial $P$, convex on $\R^d$, of degree
at most $L_d m$, such that
\begin{equation}\label{eq:Kfunctional}
 \norm{f-P}_{\B}\le c\bigl(\norm{f-h}_{\B}
 +m^{-2}\norm{D^2h}_{\B,\op}\bigr).
\end{equation}
\end{lemma}
\begin{proof}
Put $\varepsilon=\norm{f-h}_{\B}$,
$M=\norm{D^2h}_{\B,\op}$, $N=2q_d m+2$, and
\[
 \ell(x)=h(0)+\nabla h(0)\cdot x,
 \qquad q=\ell+\T_m(f-\ell).
\]
Lemma~\ref{lem:operator} shows that $q\in\mathcal P_N$ is convex on
$\B$ and
\begin{equation}\label{eq:localerror}
 \norm{f-q}_{\B}\le2\varepsilon+cm^{-2}M.
\end{equation}
The first identity in \eqref{eq:moments} gives
\[
 \T_m\ell(x)=\ell(0)+(1-\delta_m)\nabla\ell\cdot x.
\]
Thus $\T_m\ell$ is affine, and
$D^2q=D^2\T_mh+D^2\T_m(f-h)$.
By \eqref{eq:hesscontract} and \eqref{eq:Cheb}, for $r\ge1$,
\[
 |D_v^2\T_mh(re)|\le MT_N(r)|v|^2.
\]
As $\norm{\T_m(f-h)}_{\B}\le\varepsilon$,
Lemma~\ref{lem:outside} yields, for $v=ae+be^\perp$,
\begin{equation}\label{eq:negative}
 D_v^2q(re)\ge-cT_N(r)
 \bigl[(M+N^4\varepsilon)a^2+(M+N^2\varepsilon)b^2\bigr].
\end{equation}
Take $l=m+1$ and $Q_l$ from Lemma~\ref{lem:radial}. Since
$N\le2q_d l$ and
\[
 D^2[Q_l(|x|^2)]
 =2Q_l'(|x|^2)I+4Q_l''(|x|^2)xx^\top,
\]
its tangential (when $d\ge2$) and radial eigenvalues, for $|x|=r\ge1$,
are at least $c^{-1}T_N(r)$ and $c^{-1}N^2T_N(r)$,
respectively. Thus, by \eqref{eq:negative},
\[
 P(x)=q(x)+A(M+N^2\varepsilon)Q_l(|x|^2)
\]
is convex outside $\B$ for a sufficiently large constant $A$ depending only on $d$. It is convex on $\B$ as well, since $q$ is convex there and
$Q_l(|x|^2)$ is globally convex. In addition,
\[
 \norm{P-q}_{\B}\le c(M+N^2\varepsilon)l^{-2}
 \le c(Mm^{-2}+\varepsilon).
\]
Together with \eqref{eq:localerror}, this proves
\eqref{eq:Kfunctional}. Finally,
$\deg P\le8q_d l+4=8q_d m+8q_d+4\le L_d m$.
\end{proof}

\section{Proof of Theorem~\ref{thm:main}}

For $m\ge4$, choose $h=h_{1/m}$ from Lemma~\ref{lem:smoothing}
and apply Lemma~\ref{lem:global}. We obtain a polynomial $P$, convex
on $\R^d$, such that
\[
 \deg P\le L_d m,\qquad
 \norm{f-P}_{\B}\le c\,\om(f,m^{-1};\B).
\]
We have
\begin{equation}\label{eq:scale}
 \om(f,kt;\B)\le k^2\om(f,t;\B),\qquad k\in\mathbb N.
\end{equation}
This follows by restricting $f$ to each chord of $\B$ and using the
standard scaling property of the univariate second modulus; see
\cite[Section~2.7]{DeVoreLorentz1993}.
For $n\ge4L_d$ take $m=\lfloor n/L_d\rfloor$. Then $m\ge4$,
$\deg P\le n$, and $m^{-1}\le2L_d/n$; hence \eqref{eq:main}
follows from \eqref{eq:scale}.

For the remaining degrees, we use the supporting-affine-function
argument from \cite[proof of Theorem~1.3, Section~2.2]{KairePrymak2025}.
We retain the short proof for completeness.

For $1\le n<4L_d$, let $\ell$ be a supporting affine function for
$f$ at the origin. If $g=f-\ell$, then $g\ge0$, $g(0)=0$, and
\[
 0\le g(x)\le g(x)+g(-x)
 =f(x)+f(-x)-2f(0)\le\om(f,1;\B).
\]
Thus
$\norm{f-\ell}_{\B}\le\om(f,1;\B)
\le n^2\om(f,n^{-1};\B)$, and the finite range of small
$n$ is absorbed into the dimension-dependent constant.\qed

\section*{AI assistance}
ChatGPT was used to assist in developing the mathematical
argument and editing the manuscript. The authors take full responsibility for the content of the work.

\end{document}